\documentclass[english,a4paper,10pt]{article}
\usepackage[verbose, a4paper]{geometry}
\usepackage{graphicx}%
\usepackage{multirow}%

\usepackage{amsthm,mathtools,amssymb,nicefrac}%
\mathtoolsset{showonlyrefs}
\usepackage{booktabs}  
\usepackage[unicode=true, bookmarks=true, bookmarksnumbered=false, bookmarksopen=false, breaklinks=false, pdfborder={0 0 0}, pdfborderstyle={}, backref=page, colorlinks=true, pdfauthor={Rajmadan Lakshmanan}, citecolor=ForestGreen, urlcolor= darkgray, linkcolor= ForestGreen]{hyperref}
\usepackage{mathrsfs}
\usepackage[round]{natbib}
\usepackage{enumitem}
\usepackage{nowidow}
\usepackage{listings}
\usepackage{upquote}  
\setlist[enumerate,1]{label=(\roman*)}
\usepackage{textcomp}
\usepackage{nowidow}
\usepackage{dsfont}
\usepackage{xcolor}
\usepackage[ruled, vlined]{algorithm2e}	
\usepackage{todonotes}
\usepackage{csquotes} 
\usepackage{caption}
\usepackage{tikz}
\usepackage{pgfplots}
\usepackage{subcaption}
\pgfplotsset{compat=1.18}
\usetikzlibrary{positioning, calc, arrows.meta}
\usepackage{xcolor}%

\definecolor{ForestGreen}{HTML}{005F50} 

\lstdefinelanguage{Julia}{
  keywords={using, function, return, if, else, elseif, end, for, while, break, continue, struct, mutable, begin, quote, try, catch, finally, module, export, import, do, true, false, abstract, type, let, local, global, const},
  sensitive=true,
  morecomment=[l]{\#}, 
  morestring=[b]",      
}

\theoremstyle{plain}

\theoremstyle{definition}

\theoremstyle{remark}
	
\numberwithin{equation}{section}   

\usepackage{microtype}
\usepackage{doi}
\usepackage{orcidlink}

\DeclareMathOperator{\E}{{\mathds E}}

\newcommand{\one}{{\mathds 1}} 		

\title{StochasticDominance.jl: A Julia Package for Higher Order Stochastic Dominance}
\author{%
  Rajmadan Lakshmanan%
		\thanks{Faculty of Mathematics, University of Technology, Chemnitz, Germany}\,
		\footnote{\orcidlink{0009-0006-3273-9063} \href{https://orcid.org/0009-0006-3273-9063}{https://orcid.org/0009-0006-3273-9063}. Contact: \protect\href{rajmadan.lakshmanan@math.tu-chemnitz.de}{rajmadan.lakshmanan@math.tu-chemnitz.de}}
	\and
	Alois Pichler%
		\footnotemark[1]\,
		\thanks{\href{https://orcid.org/0000-0001-8876-2429}{\orcidlink{0000-0001-8876-2429} https://orcid.org/0000-0001-8876-2429}	
		}
}

\begin{document}
	\maketitle
 
\section*{Summary}
Stochastic dominance plays a key role in decision-making under uncertainty and quantitative finance.  
It compares random variables using their distribution functions.  
This concept helps to evaluate whether one investment, policy, or strategy is better than others in uncertain conditions.  
By using cumulative distribution functions, stochastic dominance allows decision-makers to make comparisons without assuming specific utility functions.  
It provides a mathematically rigorous method often used in optimization to maximize returns or minimize risk.
Its precision and reliability make it a powerful tool for analyzing complex probabilistic systems.

Despite being a crucial tool in decision-making under uncertainty and quantitative finance, (higher-order) stochastic dominance involves infinitely many constraints, making it \emph{computationally intractable} in practice.  
Our recent research \citet{lakshmanan2025HSDPrePrint} addresses this challenge by theoretically reducing the infinite constraints to a finite number.  
Building on this theoretical foundation, we have developed algorithms to maximize expected returns and minimize risk by satisfying the (higher-order) infinitely many stochastic dominance constraints.  
The paper aims to establish the mathematical completeness of this reduction and improve accessibility in research domains.  
However, no concrete, user-friendly implementation of (higher-order) stochastic dominance has been developed.
Additionally, the existing prominent theoretical algorithms only discuss stochastic orders \emph{two} and \emph{three}, but not higher orders.
Moreover, both, the discussion and implementation of non-integer orders, are absent.

To address this gap, we present \emph{StochasticDominance.jl}, an open-source Julia package tailored for verification and optimization under higher-order stochastic dominance constraints.  
Designed as a black-box solution, it enables users to input their data and obtain results effortlessly, without requiring extensive technical knowledge.

There is an extensive body of research on the applicability of stochastic dominance across various domains.  
Numerous scholars, including \citet{RuszOgryczak}, \citet{dentcheva2003optimization}, \citet{kuosmanen2004efficient}, \citet{ThirdKopa}, \citet{kopa2023multistage}, \citet{maggioni2016bounds, maggioni2019guaranteed}, and \citet{consigli2023asset}, have rigorously explored its significance.  
Their work highlights the critical role of stochastic dominance in reducing computational complexity and establishing it as a powerful tool for solving complex optimization problems.  
For a comprehensive recent textbook overview, see \citet{dentcheva2024risk}.


\section*{Main features of the package}
The StochasticDominance.jl package offers robust functions for verifying higher-order stochastic dominance constraints between two random variables.
Additionally, it provides optimization capabilities to construct an optimal portfolio that satisfies and confirms higher-order stochastic dominance constraints. Detailed explanations and usage guidelines are available in the tutorials section.

\paragraph{Technical highlights.}
StochasticDominance.jl supports two primary objective functions by maximizing expected returns and minimizing higher-order risk measures to achieve the optimal asset allocation while satisfying higher-order stochastic dominance constraints.
The package's optimization framework is built around Newton’s method, which efficiently handles the non-linear constraints.  
To enhance efficiency, we first employ Particle Swarm Optimization (PSO), which approximates the solution over a set number of iterations.  
In our previous work, \citet{lakshmanan2025HSDPrePrint} initially impose two fixed higher-order stochastic dominance constraints and dynamically introduce additional constraints to ensure dominance.
To simplify the process and align with a black-box approach, this package fixes these constraints with additional theoretical backing, eliminating the need for dynamic adjustments.
The solutions obtained from PSO serve as the initialization for Newton’s method, improving both convergence and accuracy.
Below, we provide a concise overview of its key functions.

\begin{enumerate}
 \item `verify$\_$dominance': This function checks whether the given benchmark asset, represented as the random variable $X$, and the weighted portfolio asset, represented as the random variable $Y$, exhibit a dominance relationship for the specified stochastic order. This means that $Y$ consistently yields preferable outcomes over $X$ in the specified stochastic order.

\item `optimize$\_$max$\_$return$\_$SD': This function determines the optimal asset allocation that maximizes expected returns for a given stochastic order (`SDorder'). Additionally, using `optimize$\_$max$\_$return$\_$SD(; plots=true)', users can generate a pie chart displaying the optimal allocation in percentages, along with the maximized expected returns and benchmark returns.
The function also includes the option `optimize$\_$max$\_$return$\_$SD(; verbose=true)', which allows users to imprint the convergence (or dominance) of the numerical method.

\item `optimize$\_$min$\_$riskreturn$\_$SD': This function determines the optimal asset allocation by minimizing higher-order risk measures for a given stochastic order (`SDorder') while also indicating whether dominance is achieved.
Additionally, using `optimize$\_$min$\_$riskreturn$\_$SD(; plots=true)', users can generate a pie chart that visualizes the optimal allocation in percentages, along with the minimizing higher-order risk measure returns.
The function also provides the option optimize$\_$min$\_$riskreturn$\_$SD(; verbose=true), allowing users to assess the convergence (or dominance) of the numerical algorithm.
\end{enumerate}


Further implementation details and various examples are available in the package’s documentation.\footnote{Cf.\ \href{https://rajmadan96.github.io/StochasticDominance.jl/dev/}{https://rajmadan96.github.io/StochasticDominance.jl/dev/}}
The next section presents a typical example of how to use StochasticDominance.jl.

\section{Example: verification and portfolio optimization}
The first example demonstrates how to use the `verify$\_$dominance' function to check stochastic dominance of a specified order.
The first step involves loading the necessary library, importing the data, and defining the parameters in Julia as follows:
\newpage
\begin{lstlisting}
  # Load the package
  julia> using StochasticDominance
  # Load the data 
  julia> Y = [3, 5, 7, 9, 11] # Random variable Y 
  julia> p_Y = [0.15, 0.25, 0.30, 0.20, 0.10] # Probabilities associated with Y
  julia> X = [2, 4, 6, 8, 10]   # Random variable X 
  julia> p_X = [0.10, 0.30, 0.30, 0.20, 0.10] # Probabilities associated with X
\end{lstlisting}
Next, define the stochastic order and execute the function as shown below
\begin{lstlisting}
  # Define stochastic order
  julia> SDorder = 2;
  # Run the function verify_dominance
  julia> verify_dominance(Y, X, SDorder;p_Y,p_X)
  "Y dominates X in stochastic order 2" 
\end{lstlisting}
This function checks whether $Y$ stochastically dominates $X$ of order `SDorder'. Enabling the `verbose=true' option provides detailed insights into the dominance verification process.
In this case `Y dominates X in stochastic order 2'. 

\subsection{Optimization: maximize expected return}
Next, we demonstrate how to find the optimal allocation that maximizes the expected return of the portfolios of interest while satisfying stochastic dominance of a given order.
This problem falls under non-linear optimization and involves infinitely many constraints. For a more comprehensive technical explanation, we recommend referring to the tutorials.\footnote{Cf.\ \href{https://rajmadan96.github.io/StochasticDominance.jl/dev/tutorial/tutorial3/}{https://rajmadan96.github.io/StochasticDominance.jl/dev/tutorial/tutorial3/}}

Now, we load the dataset and define the necessary parameters. We use the dataset from \citet{fama2023production}, which is a robust dataset reflecting real-world cases.
\begin{lstlisting}
  # Load the data (Fama and French 2023)
  julia> using Dates, DataFrames
  julia> data = DataFrame(
    Date = Date.([Date("2024-07-01"), Date("2024-07-02"), Date("2024-07-03"), Date("2024-07-05"), Date("2024-07-08"), Date("2024-07-09"), Date("2024-07-10"), Date("2024-07-11"), Date("2024-07-12"), Date("2024-07-15"), Date("2024-07-16"), Date("2024-07-17"), Date("2024-07-18"), Date("2024-07-19"), Date("2024-07-22"), Date("2024-07-23"), Date("2024-07-24"), Date("2024-07-25"), Date("2024-07-26"), Date("2024-07-29"), Date("2024-07-30"), Date("2024-07-31")]),
    Asset_1 = [-1.01, -2.50, -0.38, -1.11, 0.44, 0.05, -0.34, 4.00, 1.76, 1.53, 2.77, 0.71, -2.24, 0.08, 0.35, 2.55, -2.21, 1.67, 0.17, -0.97, -0.11, 0.24],
    Asset_2 = [-0.72, -0.22, 0.27, 0.15, -0.22, -1.49, 0.79, 1.60, 1.04, 0.02, 1.28, 0.69, -1.24, -1.51, 0.29, -0.09, 2.88, -0.21, 1.26, -0.79, 0.79, 1.89],
    Asset_3 = [1.10, -0.86, -0.21, -0.33, -0.64, 1.52, 1.31, 2.53, -0.06, -3.42, 2.59, -1.44, -1.74, -0.23, -0.54, 0.08, -1.61, 1.05, 2.29, -0.48, 0.98, 0.72],
    Asset_4 = [-1.80, -0.09, 0.41, 2.10, 0.59, -2.59, 2.75, 0.73, 0.41, -2.50, 0.40, 1.12, -1.72, -2.67, -0.42, -0.39, -2.51, 0.00, 0.93, -0.48, -1.65, -1.14],
    Asset_5 = [-0.65, 0.64, 0.41, -1.61, 0.33, -0.26, 1.93, 3.72, 2.17, -0.63, 1.83, 0.94, -1.78, -1.81, -0.25, 1.14, -0.41, 2.21, 2.13, -1.13, 0.07, -1.23])
  julia> xi =  Matrix(select(data, Not(:Date)))'  # Define Portfolio matrix (d=5 assets and n=22 scenarios)
  julia> d, n = size(ξ) 
  julia> tau = fill(1/d,d)  # equally weights
  julia> xi_0 = vec(tau'*xi) # Define Benchmark  
  julia> p_xi = fill(1/n,n)  # Probability vectors for x'*xi, where x is portfolio weights 
  julia> p_xi_0 = fill(1/n,n) # Probability vectors for xi_0 
\end{lstlisting}
Use the following function to compute the optimal allocation (objective: maximize expected return):
\begin{lstlisting} 
  # Define stochastic order
  julia> SDorder = 4;
  julia> x_opt, t_opt = optimize_max_return_SD(xi, xi_0,SDorder;p_xi, p_xi_0,plot=true,verbose=true) # Run the optimization
  "Simplex Constraints residuals: 5.932946379516579e-6
  Stochastic Dominance Constraints residuals: 0.0"
\end{lstlisting}
From x$\_$opt, we obtain the optimal asset allocation.
From t$\_$opt, we obtain the optimal $t$ ensuring stochastic dominance of the specified order.
For further technical details, refer to \cite{lakshmanan2025HSDPrePrint}.
The verbose option provides information about convergence status, simplex constraint residuals, and stochastic dominance constraint residuals for the given order.
Enabling plot=true generates graphical representations summarizing key insights concisely.
\begin{figure}[h]
  \centering
  \includegraphics[width=0.8\textwidth]{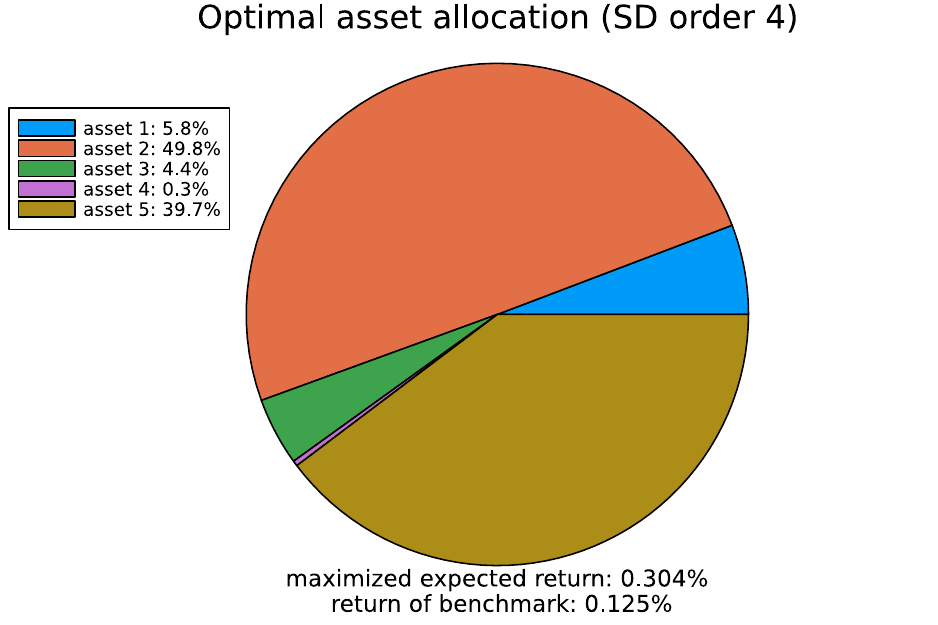}
  \caption{Optimal asset allocation (SD order 4.0). The chart illustrates the proportion of different assets in an optimized portfolio, where asset~2 holds the highest allocation (49.8\%), while asset~4 has minimal allocation. The maximized expected return is 0.304\%, compared to a benchmark return of 0.125\%.}
  \label{fig:sd_order_assets}
\end{figure}

If no portfolio allocation satisfies the stochastic dominance constraint of the given order, the algorithm stops and provides the necessary information.
\subsection{Optimization: minimizing higher-order risk measures}

Next, we demonstrate how to determine the optimal allocation that minimizes higher-order risk measures while ensuring stochastic dominance of a given order for the portfolios of interest.
This problem belongs to the class of non-linear optimization. For a detailed definition and technical explanation, we recommend referring to the tutorials.\footnote{Cf.\ \href{https://rajmadan96.github.io/StochasticDominance.jl/dev/tutorial/tutorial4/}{https://rajmadan96.github.io/StochasticDominance.jl/dev/tutorial/tutorial4/}}

We use the same dataset from the previous section and define the necessary parameters:
\begin{lstlisting}
  julia> beta=0.5 #risk parameter
  julia> r=2.0 #order of risk measure
\end{lstlisting}
Use the following function to compute the optimal allocation (objective: minimizing higher-order risk measure):
\begin{lstlisting}
  # Define stochastic order
  julia> SDorder = 4.7;
  julia> x_opt, q_opt, t_opt = optimize_min_riskreturn_SD(xi, xi_0,SDorder;,p_xi, p_xi_0,beta,r,plot=true) # Run the optimization
\end{lstlisting}
From x$\_$opt, we obtain the optimal asset allocation.
From q$\_$opt, we derive the optimal parameter that satisfies the risk measure for the given portfolio.
It is important to note that the objective itself is an optimization problem. However, our algorithm is designed to compute both simultaneously in a single execution.
Additionally, the algorithm supports non-integer stochastic dominance orders.
This feature is not restricted to this single setup but is available across all functions.
\begin{figure}[h]
  \centering
  \includegraphics[width=0.8\textwidth]{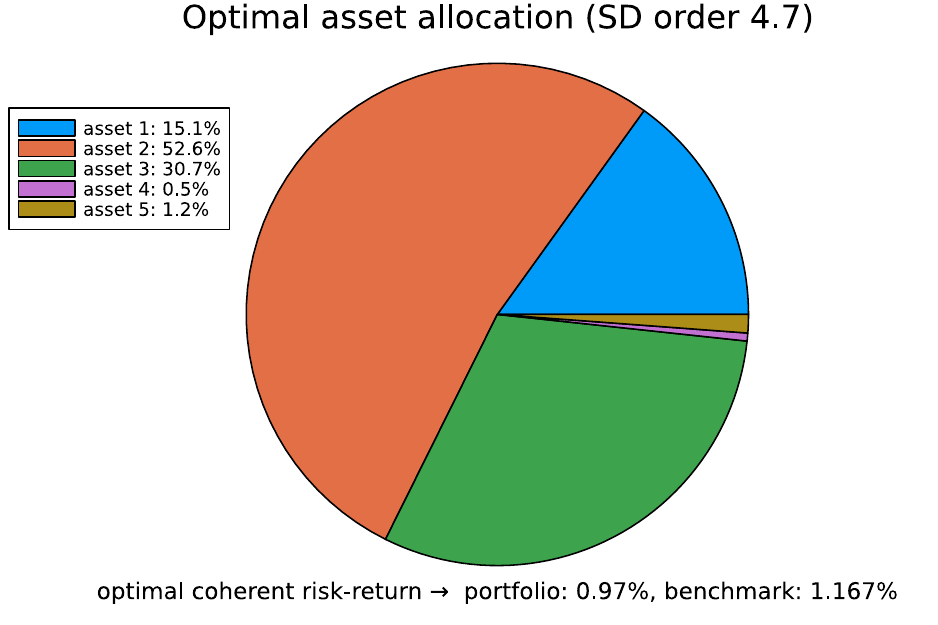}
  \caption{Optimal asset allocation (SD order 4.7). The pie chart illustrates the distribution of assets in the optimized portfolio, with asset~2 having the highest allocation (52.6\%) and asset~4 having the smallest allocation (0.5\%). The optimal coherent risk-return balance results in a portfolio return of 0.97\%, compared to a benchmark return of 1.167\%.}
  \label{fig:sd_order_assets2}
\end{figure}

\bibliographystyle{abbrvnat}
\bibliography{LiteraturAlois.bib,LiteraturRaj.bib}
\end{document}